\documentclass[10pt,a4paper]{article}
\usepackage{amsmath,amssymb,amsthm}
\usepackage[T1]{fontenc}
\usepackage[latin1]{inputenc}

\theoremstyle{plain}

\renewcommand{\H}{\mathcal{H}}
\newcommand{\F}{\mathcal{F}}
\renewcommand{\P}{\mathbb{P}}

\newcommand{\E}{\mathbb{E}}
\newcommand{\I}{1{\hskip -2.5 pt}\hbox{I}}

\begin{document}

\title{Subordination of Predictable Compensators}
\author{Henry Chiu\footnote{chiu@math.hu-berlin.de; Institut für Mathematik, Humboldt-Universität zu Berlin}}
\date{August 9, 2015}
\maketitle

\begin{abstract}
	We consider general subordination and obtain the formula of the \emph{subordinated} predictable compensator. An example of application is given.
\end{abstract}

\begin{center}\section*{\normalsize{Introduction}}\end{center}

	The idea of subordination (i.e. \hspace{-0.03cm}obtaining a new stochastic process by a random time change) was first introduced by Salomon Bochner in 1949 and is widely applied in the modelling of random phenomena such as stock price movements (e.g. the old Wall Street adage that "it takes volume $Z$ to move prices $X_Z$"). In many applications, the \emph{subordinated} process $X_Z$ is discontinuous.
	
	A central tool in the study of discontinuous process is the predictable compensator 
that arose from the general theory of stochastic processes [1]. The predictable compensator, which can be seen as a generalisation of the Lévy measure, gives a tractable description of the jump structure of a general stochastic process. It is an indispensable tool in many applications, for example, when performing an equivalent change of measure, an important operation in financial mathematics. 
	
	For a general time changed Markov process, the formula of the associated predicable compensator is not known. The purpose of this paper is to obtain such a formula.\newline
	
\begin{center}\textbf{Results}\end{center}

	It is widely known that when the time of a Lévy process $X$ is changed by an increasing Lévy process $Z$ independent of $X$, the \emph{subordinated} process $X_Z$ is a Lévy process and the \emph{subordinated} predictable compensator $(\mu^{X_{Z}})^{\P}$ of the random jump measure of $X_Z$ can be obtained by [2,Thm 30.1]: \begin{eqnarray}(\mu^{X_{Z}})^{\P}(dt,dy)=\gamma(\mu^{X})^{\P}(dt,dy)+\int_{\mathbb{R}_{+}}P^{X}_z(dy)(\mu^{Z})^{\P}(dt,dz),\end{eqnarray}where $P^{X}_t(dy)$ is the distribution of $X_t$ and $\gamma t=Z_t-\sum_{s\leq t}{\Delta{Z_s}}$. Extension of (1) to the case where $Z$ is an additive subordinator has been considered in [5,Prop.1].\newline
	
	When $X$ is replaced by, for example, more general diffusion process, (1) will no longer be applicable (the \emph{subordinated} predictable compensator shall no longer be deterministic). We extend (1) to the case where $X$ is a quasi left-continuous [1,Def.I.2.25] Markov process and $Z$ is an increasing process and give an example of application.\newline

\begin{center}\textbf{Definitions and Framework}\end{center}

	Let $X$ and $Z$ be two independent real-valued càdlàg processes defined on a complete probability space $(\Omega,\mathbb{F},\P)$ and $Z$ be increasing (i.e. non-decreasing). Denote $X_Z$ for the process obtained by a time-change of $X$ by $Z$. Let $\F$ be a right-continuous filtration in $\mathbb{F}$ such that $X_Z$ is $\F$-adapted, a non-negative random measure $(\mu^{X_Z})^{\P}$ on $\mathcal{B}(\mathbb{R}_{+}\times\mathbb{R})$ is called the $\F$-predictable compensator of the random jump measure of $X_Z$ [1,Thm.II.1.8.(i)] \& [1,Thm.I.2.2.(i)] if for all $\F$-stopping times $T$ and $B\in\mathcal{B}(\mathbb{R}\backslash\{0\})$\begin{eqnarray*}\E\sum_{t\leq T}\I_{B}(\Delta{(X_Z)}_t)=\E\int_0^{T}\int_{B}(\mu^{X_Z})^{\P}(\omega,dt,dy)\end{eqnarray*}and that the integral process $\int_0^{t\wedge T}\int_{B}(\mu^{X_Z})^{\P}(\omega,ds,dy)$ is $\F$-predictable.\\

	Denote (and respectively for $Z$ and $X_Z$) $\F^{X}:=(\F^{X}_t)_{t\geq 0}$ for the right-continuous and completed canonical filtration of $X$, $\F^{X}_{t-}:=\bigvee_{s<t}\F^{X}_s$, $\mathcal{P}^{X}$ for the $\F^{X}$-predictable $\sigma$-algebra on $\Omega\times\mathbb{R}^{+}$ generated by all left-continuous and $\F^{X}$-adapted processes and $\P^{X}$ for the probability measure $\P$ reduced to $\F^{X}_{\infty}$. Denote (and respectively for $Z$) $(\mu^{X})^{\P}$ for the $\F^{X}$-predictable compensator of random jump measure of $X$. For a $\mathcal{P}^{X}\otimes\mathcal{B}(\mathbb{R})$-measurable function $W$, we write $W\ast(\mu^{X})^{\P}$ for the integral process $W\ast(\mu^{X})^{\P}_t:=\int_0^{t}\int_{\mathbb{R}}W(\mu^{X})^{\P}(\omega,dt,dy)$\\

	Let $\mathbb{D}$ denote the space of real-valued càdlàg functions $t\mapsto u(t)$ on $\mathbb{R}^{+}$ and $\mathcal{D}_t$ denote the right-continuous smallest $\sigma$-algebra on $\mathbb{D}_t$ generated by the collection of maps $\left\{u\mapsto u(s)\right\}_{s\leq t}$, $\mathcal{D}_{t-}:=\bigvee_{s<t}\mathcal{D}_s$  and $\mathcal{D}:=\bigvee_{t\geq 0}\mathcal{D}_t$. The product space $(\mathbb{D}\times\mathbb{D},\mathcal{D}\otimes\mathcal{D})$ shall be denoted by $(\mathbb{D}^{\times 2},\mathcal{D}^{\otimes 2})$. For $u,v\in\mathbb{D}$, the map $(u,v)\mapsto(u\circ v)$ is $(\mathcal{D}^{\otimes 2},\mathcal{D})$-measurable. We shall write $(u\circ v)(t)$ for $u(v(t))$ and $(u\circ v)(t-)=u(v(t-)-)$ for $\lim_{s\uparrow t}(u\circ v)(s)$ with monotonic increasing $v$ and write\begin{eqnarray}\E^{X}F\left(X_{\cdot},Z_{\cdot}\right):=\int_\Omega F\left(X_{\cdot}(\omega),Z_{\cdot}(\widetilde{\omega})\right)\P^{X}(d\omega)\end{eqnarray}for all $(\mathcal{D}^{\otimes 2},\mathcal{B}(\mathbb{R}^{+}))$-measurable function $F$, where $\widetilde{\omega}$ is any element in $\Omega$ that is being held fixed. We define the time-changed process $X_Z$ by
\begin{eqnarray*}(X_Z)_t(\omega):=X_{Z_t(\omega)}(\omega)\end{eqnarray*}and the filtration $\F$ by \begin{eqnarray}\F_t:=(\F^{X_Z}_{t}\bigvee\F^{Z}_t)^{\P}\end{eqnarray} and $\H$ by $\H_t:=(\F^{X}_{\infty}\bigvee\F^{Z}_t)^{\P}$. Denote $\mathcal{P}$ and $\mathcal{Q}$ for the $\F$ (resp. $\H$)-predictable $\sigma$-algebra on $\Omega\times\mathbb{R}_+$. We observe $\mathcal{P}\subset\mathcal{Q}$. We shall also denote $\mathcal{P}^{X}_{Z(\widetilde{\omega})}$ for the predictable $\sigma$-algebra on $\Omega\times\mathbb{R}_{+}$ taken with respect to the filtration \begin{eqnarray}\F^{X}_{Z(\widetilde{\omega})}:=(\F^{X}_{Z_t(\widetilde{\omega})})_{t\geq 0}\end{eqnarray} for every $\widetilde{\omega}\in\Omega$ \emph{held fixed} and call a set $N\in\Omega\times\mathbb{R}_{+}$ $\P$-evanescent if $\{\omega\in\Omega: \exists t\in\mathbb{R}_+,(\omega,t)\in N\}$ is $\P$-null. If $X$ is a Markov process, we write \begin{eqnarray}P^{X}_{t}(x,s,dy):=\P(X_{s+t}\in dy|X_s=x).\end{eqnarray}

\noindent {\bf Proposition 1\it} $\;$\\
Let $A(\omega)\geq 0$ be $\H_{t}$ (resp. $\H_{t-}$)-measurable, then there exists  a $\mathcal{D}^{\otimes 2}$-measurable $H(u,v)\geq 0$ such that $A(\omega)=H(X_{\cdot}(\omega),Z_{\cdot}(\omega))$ $\P$-a.s. and
\begin{eqnarray}
\omega\longmapsto H(X_{\cdot}(\widetilde{\omega}),Z_{\cdot}(\omega))
\end{eqnarray} is $\F^{Z}_{t}$ (resp. $\F^{Z}_{t-}$)-measurable for $\P$-a.s. $\widetilde{\omega}\in\Omega$ held fixed.
If in addition, $A$ is $\F_{t}$ (resp. $\F_{t-}$)-measurable, then 
\begin{eqnarray}
\omega\longmapsto H(X_{\cdot}(\omega),Z_{\cdot}(\widetilde{\omega}))
\end{eqnarray} is $\F^{X}_{Z_{t}(\widetilde{\omega})}$ (resp. $\F^{X}_{Z_{t-}(\widetilde{\omega})-}$)-measurable for $\P$-a.s. $\widetilde{\omega}\in\Omega$ held fixed.\\
\begin{proof}
$\Pi^{1}_t(\omega):=(X_{\cdot}(\omega),Z_{\cdot\wedge t}(\omega))$, then
$\Pi^{1}_t$ is a random variable defined on $(\Omega,\H_t,\P)$ taking values in $(\mathbb{D}^{\times 2},\mathcal{D}^{\otimes 2})$ and one sees $\H_t=((\Pi^{1}_t)^{-1}\mathcal{D}^{\otimes 2})^{\P}$ by the construction of $\H_t$. If we denote $Z_{\cdot\wedge t-}$ for the map $(s\mapsto Z_{s\wedge t-})\in\mathbb{D}$, then $\H_{t-}=((\Pi^{1}_{t-})^{-1}\mathcal{D}^{\otimes 2})^{\P}$ and hence if $A(\omega)=\sum a_i\I_{A_i}(\omega)$ for $A_i\in\H_{t}$ (resp. $\H_{t-}$) then $A(\omega)=\sum a_i\I_{B_i}(\Pi^{1}_t(\omega))$ (resp. $\I_{B_i}(\Pi^{1}_{t-}(\omega))$) $\P$-a.s. for some $B_i\in\mathcal{D}^2$. The first claim holds on simple $A$.

If in addition, $A_i\in\F_{t}$ (resp. $\F_{t-}$), we define a $\mathcal{D}^{\otimes 2}$-measurable map $\Pi^{2}_t(u,v):=(u\circ v, v)(\cdot\wedge t)$ and $\Pi(\omega)_t:=(\Pi^{2}_t\circ\Pi^{1}_t)(\omega)=(X_{Z_{\cdot\wedge t}(\omega)}(\omega),Z_{\cdot\wedge t}(\omega))$. Observe also $(s\mapsto (u\circ v)(s\wedge t-))\in\mathbb{D}$, one sees $\F_{t}=(\Pi_t^{-1}\mathcal{D})^{\P}$ and $\F_{t-}=(\Pi_{t-}^{-1}\mathcal{D})^{\P}$ hence, $A(\omega)=\sum a_i\I_{B_i}(\Pi_t(\omega))$ (resp. $\I_{B_i}(\Pi_{t-}(\omega))$) $\P$-a.s. for some $B_i\in\mathcal{D}^2$. Since
$\I_{B_i}(X_{Z_{\cdot\wedge t}(\omega)}(\omega),Z_{\cdot\wedge t}(\omega))=\I_{B_i}(\Pi^{2}_t(X_{\cdot}(\omega),Z_{\cdot\wedge t}(\omega)))$
and that the path $s\mapsto X_{Z_{s\wedge t-}(\omega)}(\omega)\equiv X_{Z_{t-}(\omega)-}(\omega)$ for $s\geq t$, we see that the second claim also holds on simple $A$.

If $(H_n)_{n\geq 1}$ and $H$ are $\mathcal{D}^{\otimes 2}$-measurable, then $H_n\circ\Pi_t^{1}\rightarrow H\circ\Pi_t^{1}$ $\P$-a.s. on $\Omega$ $\Leftrightarrow$ $H_n\rightarrow H\quad\P\circ(\Pi_t^{1})^{-1}$-a.s. on $\mathbb{D}^{\times 2}$. By a monotone class argument, the claims follow.
\end{proof}
\quad

\noindent {\bf Proposition 2\it} $\;$\\
Let $A(\omega,t)\geq 0$ be $\mathcal{Q}$-measurable, then there exists a $\mathcal{D}^{\otimes 2}\otimes\mathbb{R}_{+}$-measurable $H((u,v),t)\geq 0$ such that 
$A(\omega,t)=H(X_{\cdot}(\omega),Z_{\cdot}(\omega),t)$ up to a $\P$-evanescent set and
\begin{eqnarray}
(\omega,t)\longmapsto H(X_{\cdot}(\widetilde{\omega}),Z_{\cdot}(\omega),t)
\end{eqnarray} is $\mathcal{P}^{Z}$-measurable for $\P$-a.s. $\widetilde{\omega}\in\Omega$ held fixed.
If in addition, $A$ is $\mathcal{P}$-measurable, then 
\begin{eqnarray}
(\omega,t)\longmapsto H(X_{\cdot}(\omega),Z_{\cdot}(\widetilde{\omega}),t)
\end{eqnarray} is $\mathcal{P}^{X}_{Z(\widetilde{\omega})}$-measurable for $\P$-a.s. $\widetilde{\omega}\in\Omega$ held fixed.

\begin{proof}
The claims clearly holds for all $\mathcal{Q}$-measurable (resp. $\mathcal{P}$-measurable) $A$ of the form $A_t=A_0\I_{\{0\}}(t)+\sum_{i\in\mathbb{N}} A_{t_i}\I_{(t_i,t_{t+1}]}(t)$ for $\H_{t_i-}$ (resp. $\F_{t_i-}$)-measurable $A_{t_i}$ as a direct consequence of Proposition 1. Observe also if $(H_n)_{n\geq 1}$ and $H$ are $\mathcal{D}^{\otimes 2}\otimes\mathbb{R}_+$-measurable then $H_n(\Pi_t^{1}(\omega),t))\rightarrow H(\Pi_t^{1}(\omega),t))$ on $\Omega\times\mathbb{R}_+$ up to a $\P$-evanescent set $\Leftrightarrow$ $H_n((u,v),t)\rightarrow H((u,v),t)$ on $\mathbb{D}^{\times 2}\times\mathbb{R}_{+}$ up to a $\P\circ(\Pi_t^{1})^{-1}$-evanescent set. By a monotone class argument, the claims follow.
\end{proof}
\quad

\noindent {\bf Theorem\it} $\;$\\
\textit{Let $X$ be a quasi left-continuous Markov process with transition kernel $P^{X}_{t}(x,s,dy)$ and $Z$ be an increasing process independent of $X$. Denote $X_Z$ the process obtained by a time-change of $X$ by $Z$ and $Z^{c}_{t}:=Z_{t}-\sum_{s\leq t}{\Delta{Z_s}}$ then $(\mu^{X_{Z}})^{\P}(\omega,dt,dy)$ is changed as follows:}\\
\begin{eqnarray}
(\mu^{X})^{\P}(\omega,dZ^{c}_t,dy)+\int_{\mathbb{R}_{+}}P^{X}_z(X_{Z_{t-}},Z_{t-},\{X_{Z_{t-}}\}+dy)(\mu^{Z})^{\P}(\omega,dt,dz).
\end{eqnarray}\quad

\begin{proof}
Let $A\times B\in\mathcal{P}\otimes\mathcal{B}(\mathbb{R}\backslash\{0\})$, $W:=\I_{A\times B}$, $I:=\{t\geq 0|\Delta{Z}_t=0\}$. Observe $\Delta{(X_Z)}_t=\Delta{(X)}_{Z_{t}}$ on $I$ and $\Delta{(X_Z)}_t=X_{Z_{t-}+\Delta{Z}_t}-X_{Z_{t-}-}$ on $I^c$ and by (2), 
we can write \begin{eqnarray*}
F(X_{\cdot},Z_{\cdot})=(W\ast\mu^{X_Z})_{\infty}=\sum_{t\in I}\I_A\I_B(\Delta{(X)}_{Z_{t}})+\sum_{t\in I^c}\I_A\I_B(X_{Z_{t-}+\Delta{Z}_t}-X_{Z_{t-}-})
.\end{eqnarray*}Let $Z^{-1}$ denote the left-continuous generalized inverse of $Z$, by (4) \& (9) put $\I_{A}=H(X_{\cdot},Z_{\cdot},t)$ then $(\omega,t)\mapsto H(X_{\cdot}(\omega),Z_{\cdot}(\widetilde{\omega}),Z^{-1}_t(\widetilde{\omega}))$ is $\mathcal{P}^{X}$-measurable for $\P$-a.s. $\widetilde{\omega}\in\Omega$ held fixed [1,Prop.I.2.4]. Together with the quasi left-continuity of $X$,  [1,Thm.II.1.8] \& [Cor.II.1.19], it follows $\E\sum_{t\in I}\I_A\I_B(\Delta{(X)}_{Z_{t}})$  (see also (2) for notation)
\begin{eqnarray*}
&=&\E^{Z}\E^{X}\sum_{t\in Z(I)}H(X_{\cdot}(\omega),Z_{\cdot}(\widetilde{\omega}),Z^{-1}_t(\widetilde{\omega}))\I_B(\Delta{X}_{t})\\
&=&\E^{Z}\E^{X}\int_{Z(I)}\int_{\mathbb{R}}H(X_{\cdot}(\omega),Z_{\cdot}(\widetilde{\omega}),Z^{-1}_t(\widetilde{\omega}))\I_B(y)(\mu^{X})^{\P}(\omega,dt,dy)\\
&=&\E^{Z}\E^{X}\int_{I}\int_{\mathbb{R}}H(X_{\cdot}(\omega),Z_{\cdot}(\widetilde{\omega}),t)\I_B(y)(\mu^{X})^{\P}(\omega,dZ_t(\widetilde{\omega}),dy)\\
&=&\E^{Z}\E^{X}\int_{\mathbb{R}_{+}}\int_{\mathbb{R}}H(X_{\cdot}(\omega),Z_{\cdot}(\widetilde{\omega}),t)\I_B(y)(\mu^{X})^{\P}(\omega,dZ^{c}_t(\widetilde{\omega}),dy)\\
&=&\E\int_{\mathbb{R}_{+}\times\mathbb{R}}W(\omega,t,y)(\mu^{X})^{\P}(\omega,dZ^{c}_t(\omega),dy)
.\end{eqnarray*}Since $X$ has no fixed times of discontinuity and that $I^{c}$ is countable and by (8), it follows $\E\sum_{t\in I^c}\I_A\I_B(X_{Z_{t-}+\Delta{Z}_t}-X_{Z_{t-}-})$

\begin{eqnarray*}
&=&\E^{X}\E^{Z}\sum_{t\in I^{c}}H\I_B(X_{Z_{t-}+\Delta{Z}_t}-X_{Z_{t-}})\\
&=&\E^{X}\E^{Z}\int_{\mathbb{R}_{+}}\int_{\mathbb{R}_{+}}H\I_B(X_{Z_{t-}+z}-X_{Z_{t-}})(\mu^{Z})^{\P}(\omega,dt,dz)\\
&=&\E^{Z}\int_{\mathbb{R}_{+}}\int_{\mathbb{R}_{+}}\E^{X}[H\I_B(X_{Z_{t-}+z}-X_{Z_{t-}})](\mu^{Z})^{\P}(\widetilde{\omega},dt,dz).\\
\end{eqnarray*}By (4), (9) and  [1,Prop.I.2.4], we see that for $\P$-a.s. $\widetilde{\omega}\in\Omega$ held fixed, the map $\omega\mapsto H(X_{\cdot}(\omega),Z_{\cdot}(\widetilde{\omega}),t)$ is $\F^{X}_{Z_{t-}(\widetilde{\omega})}$-measurable for all $t\geq 0$. Together with the Markov property of $X$ we have $\E^{X}[H\I_B(X_{Z_{t-}+z}-X_{Z_{t-}})]$

\begin{eqnarray*}
&=&\E^{X}\left[H\E^{X}[\I_B(X_{Z_{t-}(\widetilde{\omega})+z}-X_{Z_{t-}(\widetilde{\omega})})|\F^{X}_{Z_{t-}(\widetilde{\omega})}]\right]\\
&=&\E^{X}[HP^{X}_{z}(X_{Z_{t-}(\widetilde{\omega})},Z_{t-}(\widetilde{\omega}),\{X_{Z_{t-}(\widetilde{\omega})}\}+B)]
\end{eqnarray*}hence $\E\sum_{t\in I^c}\I_A\I_B(X_{Z_{t-}+\Delta{Z}_t}-X_{Z_{t-}-})$
\begin{eqnarray*}
&=&\E^{Z}\E^{X}\int_{\mathbb{R}_{+}\times\mathbb{R}_{+}}H\int_{B}P^{X}_{z}(X_{Z_{t-}},Z_{t-},\{X_{Z_{t-}}\}+dy)(\mu^{Z})^{\P}(\omega,dt,dz)\\
&=&\E\int_{\mathbb{R}_{+}\times\mathbb{R}}W\int_{\mathbb{R}_{+}}P^{X}_{z}(X_{Z_{t-}},Z_{t-},\{X_{Z_{t-}}\}+dy)(\mu^{Z})^{\P}(\omega,dt,dz).
\end{eqnarray*}Define $v(\omega,dt,dy):=$

\begin{eqnarray*}
(\mu^{X})^{\P}(\omega,dZ^{c}_{t},dy)+\int_{\mathbb{R}_{+}}P^{X}_{z}(X_{Z_{t-}},Z_{t-},\{X_{Z_{t-}}\}+dy)(\mu^{Z})^{\P}(\omega,dt,dz)
\end{eqnarray*} then $\E(W\ast\mu^{X_Z})_{\infty}=\E(W\ast v)_{\infty}$. It is clear that $v(\omega,dt,dy)$ defines a non-negative random measure on $\mathbb{R}_{+}\times\mathbb{R}$ and that $(W\ast v)_t$ is $\F$-predictable (3). If $T$ is a $\F$-stopping time, put $A:=\{(\omega,t):0\leq t\leq T(\omega)\}\in\mathcal{P}$, (10) follows.\end{proof}

\noindent {\bf Example \it} $\;$
We calculate the compensator $(\mu^{X_Z})^{\P}$ of the random jump measure of $X_Z$ with $X$ and $Z$ taken to be, respectively, a skew Brownian motion (diffusion process) and a tempered stable subordinator independent of $X$. The compensator of the random jump measure of $Z$ is 
\begin{eqnarray}
(\mu^{Z})^{\P}(dt,dz)=dt\frac{c}{z^{1+\alpha}}e^{-\lambda z}\I_{\{z>0\}}(dz)
\end{eqnarray}for $c, \lambda>0$ and $\alpha\in[0,1)$. The case $\alpha=0$ corresponds to a Gamma subordinator. By [3,(17)], the transition function of a skew Brownian motion $X$ can be written as\begin{eqnarray}
P^{X}_t(x,dy)&=&\frac{1}{\sqrt{2\pi t}}\left(e^{-\frac{(|y-x|)^2}{2t}}+\frac{\beta}{\text{sgn}(y)}e^{-\frac{(|y|+|x|)^2}{2t}}\right)dy
\end{eqnarray}for $\beta\in[-1,1]$. The case $\beta=0$ corresponds to the standard Brownian motion. Using the modified Bessel function $K_{v}(x)$ for the integral representation
\begin{eqnarray}
\int_{0}^{\infty}\frac{1}{z^{1+v}}e^{-\frac{a^2 z}{2}-\frac{b^2}{2z}}dz=2(\frac{a}{b})^{v}K_{v}(ab)
,\end{eqnarray}the compensator formula (10) and $\phi(X_{Z_{t-}}(\omega),y):=|X_{Z_{t-}}(\omega)|+|X_{Z_{t-}}(\omega)+y|$, we obtain\begin{eqnarray*}
(\mu^{X_Z})^{\P}(\omega,dt,dy)&=&\frac{2c}{\sqrt{2\pi}} \left(\frac{\sqrt{2\lambda}}{|y|}\right)^{1/2+\alpha}K_{1/2+\alpha}\left(\sqrt{2\lambda}|y|\right)dtdy\\
&+&\frac{\beta 2c}{\sqrt{2\pi}\text{sgn}(X_{Z_{t-}}+y)}\left(\frac{\sqrt{2\lambda}}{\phi
(X_{Z_{t-}},y)}\right)^{1/2+\alpha}\\
&\times&K_{1/2+\alpha}\left(\sqrt{2\lambda}\phi
(X_{Z_{t-}},y)\right)dtdy
\end{eqnarray*}and for the Gamma case $\alpha=0$,\\
\begin{eqnarray}
(\mu^{X_Z})^{\P}(\omega,dt,dy)=\left(\frac{ce^{-\sqrt{2\lambda}|y|}}{|y|}+\frac{ \beta}{\text{sgn}(X_{Z_{t-}}+y)}\frac{ce^{-\sqrt{2\lambda}\phi(X_{Z_{t-}},y)}}{\phi(X_{Z_{t-}},y)}\right)dtdy.\quad\end{eqnarray}\\
	We see that $(\mu^{X_Z})^{\P}$ is deterministic and time-independent if and only if $\beta=0$, in this case $X_Z$ is a time-changed Brownian motion. If in addition, $Z$ is a Gamma process (i.e. $\alpha=\beta=0$) then $X_Z$ is a Variance Gamma process [4] with Lévy measure $v(dy)=\frac{ce^{-\sqrt{2\lambda}|y|}}{|y|}dy$ and (14) reduces to\\
\begin{eqnarray}
(\mu^{X_Z})^{\P}(\omega,dt,dy)=dtv(dy).
\end{eqnarray}

\bibliographystyle{plain}
\bibliography{subordination_revised}

\end{document}